\documentclass[12pt]{article}
\usepackage{amssymb}

\def\beq{\begin{equation}}
\def\eeq{\end{equation}}
\def\bea{\begin{eqnarray}}
\def\eea{\end{eqnarray}}
\def\ar{\begin{array}}
\def\ear{\end{array}}

\def\nn{\nonumber}

\def\ga{\gamma}
\def\Ga{\Gamma}
\def\si{\sigma}

\def\al{\alpha}

\def\rd{{\rm d}}
\def\ri{{\rm i}}

\def\al{\alpha}

\def\ga{\gamma}

\def\si{\sigma}

\def\rd{{\rm d}}
\def\ri{{\rm i}}

\begin{document}
\begin{center}
{\bf \Large A Tail Sensitive Test\\[2mm] for Cumulative Distribution Functions}\\[1.1cm]

{\bf Krzysztof A. Meissner}

\vspace{0.3cm}
{{\it Faculty of Physics,
University of Warsaw,
Ho\.za 69, Warsaw, Poland\\
}}

\vspace{0.7cm}

\begin{minipage}[t]{11cm}
{\footnotesize
{We propose a simple way of testing whether a given set of observations can come from a given theoretical cumulative distribution. In the test more weight is attached to the tails of the distribution than in the usual Kolmogorov or Smirnov tests. The respective probability distribution is derived.}
}
\end{minipage}
\end{center}

\section{Introduction}
\noindent
In mathematical statistics it is extremely important to test whether a given sample can come from a theoretical cumulative distribution function (CDF) and to provide quantitative measures for this hypothesis. The most celebrated tests of this type are the Kolmogorov \cite{K} and Smirnov \cite{S} tests based on the supremum of the absolute distance between the observational and theoretical or two observational CDFs. For example in the Smirnov test if we have two distributions with $N_A$ and $N_B$ number of points respectively and the observational CDFs $F_A(x)$ and $F_B(x)$ then with the assumption that they come from the same CDF the probability for the weighted supremum $D$ of the absolute distance between the two, $D:=\sup_x |F_A(x)-F_B(x)|/\sqrt{N_A^{-1}+N_B^{-1}}$,  is given asymptotically by
\beq
P\left(D>\lambda\right)=2\sum_{n=1}^\infty (-1)^{n-1}\,e^{-2 n^2\lambda^2}=1-\frac{\sqrt{2\pi}}{\lambda}
\sum_{n=1}^\infty e^{-(2n-1)^2\pi^2/(8\lambda^2)}
\eeq
When $N_B\to\infty$, i.e the second distribution is treated as a theoretical one, we recover the Kolmogorov test.
Practical usage of these tests shows however that they are less sensitive to the tails of the distribution than to the bulk. In practice it is often the case that we observe exceptionally high or low values which seem to indicate that they come from a different distribution but the difference does not translate into a significant change of the confidence level in the Kolmogorov or Smirnov tests. There are tests that are more sensitive to the tails of the distribution but they usually assume normality of the distribution. We would like to propose a test which emphasizes the role of the tails of the distribution (separately for high and low ends), is independent of the normality assumption, simple to apply and has explicitly calculable statistical properties. The motivation for the introduction of the test comes from the question whether there are statistically significant patterns on the maps of the background radiation as measured by the WMAP satelite -- the description of the problem and the application of the test will be presented in a separate publication \cite{MNR}.

\section{The proposal}
Assume that the theoretical CDF is given by $F(x)$ and the ordered sample (with $n$ entries) gives the experimental CDF $F_n(x)$ (which is a step function increasing by $\frac{1}{n}$ at the points $x_i$).  The proposed test (right) is based on the following quantity
\beq
A^R_{a,n}=-\int\left(1-F_n(x)\right)\rd\ln\left(1-(F(x))^a\right)=
-\frac{a}{n}\sum_{i=1}^n\ln\left(1-(F(x_i))^a\right)
\eeq
(and analogous left test $A^L_{a,n}$ with $F\to (1-F)$) where $a$ is a positive real number. It is clear that the right test gives more weight to values of $F$ close to 1 while the left test to values close to 0. With increasing $a$ we increase the relative weight of tails of the distribution (right for $A^R$ and left for $A^L$).

The test can also be used when we use some ``coarse-graining'':  we may group the total number of observations $N$ into $n$ ($1\ll n\ll N$) bins with positions $x_i$ and $d_i$ points in the $i$th bin and use the formula
\beq
A=-\frac{a}{N}\sum_{i=1}^n d_i\ln\left(1-(F(x_i))^a\right)
\eeq

\section{Properties of the test}

To derive the distribution function for $A^R_{a,n}$ (the same for $A^L$) we first calculate (we denote by $A$ either $A^R$ or $A^L$)
\beq
\langle e^{\ri t A}\rangle=\left(\int\limits_0^1\rd z\left(1-z^a\right)^{-\ri t a/n}\right)^n
\eeq
It is straightforward to calculate this quantity with the result
\beq
\langle e^{\ri t A}\rangle=\left(\frac{\Gamma\left(1+\frac1{a}\right)
\Gamma\left(1-\frac{\ri t a}{n}\right)}{\Gamma\left(1-\frac{\ri t a}{n}+\frac{1}{a}\right)}\right)^n
\label{genfun}
\eeq
Defining cumulants of the distribution as
\beq
\ln\langle e^{\ri t A}\rangle=\sum_{k=1}^{\infty}\frac{(\ri t)^k\sigma_k(A)}{k!}
\label{momdef}
\eeq
so that for example
\bea
\si_1(A)&=&\langle A\rangle,\nn\\
\si_2(A)&=&\langle A^2\rangle-\langle A\rangle^2\nn\\
\si_3(A)&=&\langle A^3\rangle-3\langle A^2\rangle\langle A\rangle+2\langle A\rangle^3\nn
\eea
and using the properties of the $\Gamma$ function we get the general expression for the cumulants $\si_k(A)$ for the distribution (\ref{genfun}):
\beq
\si_k(A)=a\,(k-1)!\left(\frac{a}{n}\right)^{k-1}\sum_{l=1}^{\infty}\left(\frac{1}{l^k}
-\frac{1}{(l+1/a)^k}\right)
\label{momres}
\eeq
The probability distribution for $A$ is given by:
\beq
g_a(n;s)=\frac{1}{2\pi}\int\limits_{-\infty}^{\infty}\rd t\ e^{-\ri s t}\langle e^{\ri t A}\rangle
\label{dist}
\eeq
The integral gives 0 if $s<0$ which is consistent with the definition of $A$ -- in the formulae below we assume henceforth $s\ge 0$. For general $n$ and $a$ it is straightforward to calculate $g_a(n;s)$ numerically with arbitrary accuracy.

In the case ($a=1$, $n$ arbitrary) it is possible to perform the integral in (\ref{dist}) analytically and we get
\beq
g_1(n;s)=\frac{n^n s^{n-1}}{(n-1)!}e^{-ns}\equiv\Gamma_n\left(s,\frac{1}{n}\right)
\eeq
i.e. the gamma distribution with $\langle s\rangle=1$ and the variance $1/n$.

Since the distribution (\ref{dist}) is given as a Fourier transform it is straightforward to give also the expression for the cumulative distribution i.e. the probability to get the value of $A$ smaller than $\sigma$:
\beq
G_a(n;\sigma):=\int\limits^\sigma_0 \rd s\, g_a(n;s)=\frac{1}{2\pi}\int\limits_{-\infty}^{\infty}\rd t \frac{1- e^{-\ri \sigma t}}{\ri t}\langle e^{\ri t A}\rangle
\eeq
It is straightforward to calculate $G_a(n;\si)$ numerically with arbitrary accuracy.

\section{The case of large $a$ and $n$}
We define
$$
\alpha:=\frac{a}{n}
$$
and below we discuss the case  $a\to\infty,\ n\to\infty$ with $\alpha$ kept finite.

To derive the limiting distribution in this case we expand $\Gamma$ functions and we get (up to $1/n$ and $1/a$ corrections)
\beq
\langle e^{\ri t A}\rangle=\exp\left(-\frac{\gamma}{\alpha}-\frac{\psi(1-\ri t \alpha)}{\alpha}\right)
\eeq
and the distribution function
\beq
g_\infty(\alpha;s)=\frac{1}{2\pi}\int\limits_{-\infty}^{\infty}\rd t\ e^{-\ri s t-\frac{\gamma}{\alpha}-\frac{1}{\alpha}\psi(1-\ri t \alpha)}
\label{distalpha}
\eeq
where
$$
\psi(1-\ri t)=\ri\frac{\rd}{\rd t}\ln\Gamma(1-\ri t)=-\gamma-
\sum_{l=1}^{\infty}\frac{\ri t}{l(l-\ri t)}
$$
The integral cannot be calculated in a closed form but we can give a closed expression for the cumulants of the distribution using (\ref{momres})
$$
\sigma_k=k!\,\zeta(k+1)\,\alpha^{k-1}
$$

It is also convenient to give directly the cumulative distribution of $g_\infty$:
\beq
G_\infty(\alpha;\sigma)=\int\limits_0^\sigma \rd s \, g_\infty(\alpha;s)=\Re\int\limits_{0}^{\infty}\rd t\,\frac{1- e^{-\ri t}}{\ri\pi t}\, \exp\left(-\frac{\gamma}{\alpha}-\frac{1}{\alpha}\psi\left(1-\frac{\ri t \alpha}{\si}\right)\right)
\label{intdistequiv}
\eeq
It is straightforward to calculate numerically this integral for arbitrary $\alpha$ and $\sigma$ with any prescribed accuracy -- for example $G_\infty(1;1)=0.439166...$, $G_\infty(1;3)=0.8390636...$, $G_\infty(1;7)=0.9898427...$ and $G_\infty(1;17)=0.999995...$.

For $\alpha\gg 1$ it is useful to separate the part that is slowly decaying for large $t$ from the rest:
\bea
G_\infty(\al;\si)&=&
\Re\int\limits_{0}^{\infty}\rd t\,\frac{1- e^{-\ri t}}{\ri\pi t}\, \exp\left(-\frac{\gamma}{\al}-\frac{1}{\al}\ln\left(\frac{t }{y}\right)+\frac{\ri\pi}{2\alpha}\right)\nn\\
&&+\Re\int\limits_{0}^{\infty}\rd t\,\frac{1- e^{-\ri t}}{\ri\pi t}\, \left[\exp\left(-\frac{\gamma}{\alpha}-\frac{1}{\alpha}\psi\left(1-\frac{\ri t}{y}\right)\right)\right.\nn\\
&&\left.\ \ \ \ \ \ \ \ \ \ \ \ \ \ \ \ -
\exp\left(-\frac{\gamma}{\alpha}-\frac{1}{\alpha}\ln\left(\frac{t }{y}\right)+\frac{\ri\pi}{2\alpha}\right)\right]
\eea
where $y:=\frac{\si}{\al}$. The first term can be easily integrated and we get
\bea
G_\infty(\al;\si)&=&
\exp\left(\frac{\ln(y)-\ga}{\al}-\ln\left(\Ga(1+1/\al)\right)\right)\nn\\
&&+\Re\int\limits_{0}^{\infty}\rd t\,\frac{1- e^{-\ri t}}{\ri\pi t}\, \left[\exp\left(-\frac{\gamma}{\alpha}-\frac{1}{\alpha}\psi\left(1-\frac{\ri t }{y}\right)\right)\right.\nn\\
&&\left.\ \ \ \ \ \ \ \ \ \ \ \ \ \ \ -
\exp\left(-\frac{\gamma}{\alpha}-\frac{1}{\alpha}\ln\left(\frac{t }{y}\right)+\frac{\ri\pi}{2\alpha}\right)\right]
\label{ginf}
\eea
The formula can be expanded in inverse powers of $\al$ and all the integrals are well behaved.  It is convenient to organize the series in a slightly different way (the first factor reflects the leading behavior):
\beq
G_\infty(\alpha;\sigma)=\left(1-e^{-y}\right)^{1/\alpha}\left(1+\sum_{k=2}^{\infty}
\frac{f_k(y)}{\alpha^k}\right)
\label{Ginfty}
\eeq
where the functions $f_k(y)$ can be read off from (\ref{ginf}). It is rather straightforward to see that $f_1(y)=0$ and after rather involved manipulations we get:
\beq
f_2(y)=\frac{y}{2} \ln\left(1- e^{-y}\right)-\frac12\sum_{l=1}^\infty \frac{e^{-ly}}{l^2}
\eeq
so it is a monotonic function from $f_2(0)=-\frac{\pi^2}{12}$ to $f_2(\infty)=0$.
The form (\ref{Ginfty}) is very useful in actual applications \cite{MNR}.

\section{Experimental CDF}

The test described in this paper requires the knowledge of the theoretical CDF but it is often the case that we have at our disposal only measurements of CDF. If we have $k$ such measurements (each with $n$ entries) then we calculate at each point the average experimental CDF $\mu(x)$ equal to the average of all measured CDFs at this point. Therefore at every point $\mu$ can take one of the values $0,\frac{1}{kn},\frac{2}{kn},\ldots,1$. It is straightforward to prove that for the theoretical CDF given by $F(x)$ the probability to measure at this point the value $\mu$ is given by
\beq
P_n^k(\mu)=\frac{(kn)!}{(kn\mu)!(kn(1-\mu))!}F^{kn\mu}(1-F)^{kn(1-\mu)}
\eeq
For large values of $kn$ it tends to the gaussian distribution
\beq
P_n^k(\mu)\approx \sqrt{\frac{1}{2\pi kn F(1-F)}}\exp\left(-\frac{kn(\mu-F)^2}{2F(1-F)}\right)
\eeq
so that $\mu$ approximates $F$ with dispersion $\sqrt{\frac{F(1-F)}{kn}}$.

\section{Conclusions}
The test proposed in the paper requires more work to check its power and usefulness and compare it to known non-parametric tests but with its simplicity and explicit probability distributions it should be useful in statistical data analysis.

\section*{Acknowledgements}
The author thanks the AEI Max Planck Institute in Potsdam for hospitality.

\end{document}